\documentstyle[12pt,amssymbols]{article}

\def\presentation{
\voffset -.50in
\hoffset -.19in
\oddsidemargin 0in \evensidemargin 0in
\marginparwidth .75in \marginparsep 7pt \topmargin 0in
\headheight 12pt \headsep .25in
\footheight 18pt \footskip .35in
\textheight 9.5in \textwidth 6.5in
\columnsep 10pt \columnseprule 0pt }

\presentation

\newcommand{\beq}{\begin{equation}}
\newcommand{\eeq}{\end{equation}}
\newcommand{\bea}{\begin{eqnarray}}
\newcommand{\bean}{\begin{eqnarray*}}
\newcommand{\eea}{\end{eqnarray}}
\newcommand{\eean}{\end{eqnarray*}}
\newcommand{\demi}{\frac{1}{2}}
\newcommand{\sumi}{\sum_{i=1}^N}
\newcommand{\sumij}{\sum_{_{\ i \neq j}^{i,j=1}}^N}
\newcommand{\la}{\lambda}
\newcommand{\si}{\sigma}

\newcommand{\q}{q}

\newcommand{\Sum}{\sum}

\newtheorem{th}{Theorem}[section]
\newtheorem{prop}{Proposition}[section]

\newtheorem{lem}{Lemma}[section]

\begin{document}

{\bf Classical dynamical $r$-matrices for Calogero-Moser systems and their
generalizations}

\vskip0.5in

{\bf J. Avan}

\vskip0.5in

{\it L.P.T.H.E. Universit\'e Paris VI (CNRS UA 280), Box 126, Tour 16,
$1^{{\rm er}}$ \'etage, 4 place Jussieu, 75252 Paris Cedex 05, France}

\section{Introduction}

The construction and study of classical (and quantum) dynamical $r$-matrices
is currently undergoing extensive development. Various examples of such
objects were recently discussed, for instance in \cite{Tsi,Hie,KuRWTsi}. 
However at
this time one lacks a general classifying scheme such as exists in the case
of constant classical $r$-matrices thanks to Belavin and Drin'feld \cite {BZ}.
A partial classification scheme has very recently been proposed
\cite{Schiff} for dynamical
$r$-matrices obeying the particular version of the dynamical Yang-Baxter 
equation \cite{FeWie} corresponding to Calogero-Moser models \cite{ABB2}.

Consideration of such structures is thus particularly relevant to the study of
Calogero-Moser \cite{CalMo, GibHe} and relativistic Ruijsenaars-Schneider
models \cite{Rui,KriZa} where they appear systematically. Their occurence and
the particular form they assume, which we shall give in detail in this
lecture, are due to the common nature of these models as hamiltonian 
reductions of free
or harmonic motions on particular symplectic manifolds: cotangent bundle of Lie algebras or Lie groups for rational or trigonometric Calogero-Moser models
\cite{OlPer}; double of Lie groups \cite{ACF1} or cotangent bundle to centrally
extended loop group \cite{GorNe} for Ruijsenaars-Schneider models. The most 
general elliptic 
potentials are in turn associated to loop groups
over elliptic curves \cite{GorNe,ACF1} and are crucial in understanding 
the algebraic resolution of these models \cite{Kri,BBKT}.

This hamiltonian reduction procedure is also an important tool in the 
explicit construction of classical $r$-matrices for $BC_n$-type systems, where 
the direct resolution of the intricate $r$-matrices equations, such as was 
done in 
\cite{AvTaSkBraSu}, becomes untractable. We shall exemplify such a construction 
in the case of trigonometric Calogero-Moser models \cite{AvBabTa}.

The plan of this lecture runs as follows. We first recall the essential results
of classical $r$-matrix theory and introduce the notations to be used throughout
it. We then describe the construction of classical $r$-matrices for 
trigonometric Calogero-Moser models using the Hamiltonian reduction procedure.
This gives us a general formula valid for all non-exceptional Lie groups.

We finally give a systematic overall picture of the classical $r$-matrix
structure of this type which were obtained by various authors for:

a) elliptic, trigonometric and rational $A_n$ Calogero-Moser models
\cite{AvTaSkBraSu}

b) elliptic, trigonometric and rational $A_n$ spin Calogero-Moser models
\cite{BiAvBa}

c) elliptic, trigonometric and rational $A_n$ Ruijsenaars-Schneider models.
\cite{BaBe,AvRol,Su}.

The case of spin Ruijsenaars-Schneider systems \cite{KriZa} was recently 
investigated in the rational case \cite {Russes} but the general Hamiltonian
structure still escapes understanding at this time.

To be complete we must indicate that two alternative approaches were recently
described. One, using the same Lax matrices but leading to dynamical 
$r$-matrices of a different type (to be commented upon later) was developed
in \cite{Bra2} (see contribution by Pr. Braden). The other one 
uses conjugated Lax matrices which allow to eliminate the dynamical dependence 
in the $r$-matrix at the cost of introducing a more complicated Lax operator.
It was developed directly in the quantum case using the formalism of 
intertwining vectors \cite {Sacha} (see contribution by Pr. Hasegawa).

\section{Preliminaries}

First of all we need
to recall  four essential features of 
the classical $r$-matrix formalism (see \cite {FT} for a textbook presentation
of the hamiltonian theory of classical integrable systems.).

We consider a generic dynamical system described by a set of coordinates 
$\{x_i \}$ and momenta $\{p_i \}, i=1 \cdots n$; a Poisson structure 
$ \{ \} $ and a Hamiltonian $h(x_i , p_j)$.
\newline

{\bf 1- Liouville theorem}

The existence of $n$ algebraically independent, globally defined, Poisson
commuting quantities such that the hamiltonian $h$ belongs to the ring
generated by this set, guarantees the existence of a canonical transformation
$(x_i, p_j) \rightarrow (I_i, \theta_j)$ linearizing the equations of motion 
\cite{Liou}.
Further assumptions on the topological structure of the phase space allow more
precise statements on the geometrical interpretation of the transformed 
variables known as action-angle variables \cite {Ar}.
\newline

{\bf 2- Lax pair formulation}

The Lax pair formulation of a dynamical system is the giving of two elements
of a Lie algebra ${\cal G}$, $L(x,p)$ and $M(x,p)$ such that the equations of
motion for $x_i$, $p_i$ are equivalent to the isospectral evolution (Lax
equation) \cite{Lax}

\beq
\label{Lax}
\frac{dL}{dt} = [L,M]
\eeq

It follows that the adjoint-invariant quantities $Tr L^n$ are time-independent.
In order to implement Liouville theorem onto this set of possible action 
variables
we need them to be Poisson-commuting. This is ensured by the classical $r$
matrix structure.
\newline

{\bf 3- The $r$-matrix structure}

Defining the decomposition of the Lax operator on a basis $\{ t_a \}$
of the Lie  algebra ${\cal G}$ as $L \equiv \sum_{a} l^a t_a $, the Poisson 
commutation of the ad-invariant $Tr L^n$ is equivalent to the existence of
an object $r_{12}(x,p) \in {\cal G}_1 \otimes {\cal G}_2$ hereafter known as a
classical $r$-matrix \cite{Skl1,JMM1,BV1}, such that:

\beq
\{ L_1 , L_2 \} \equiv \sum_{a,b} \{l^a , l^b\} t_a \otimes t_b =
[r_{12}, L_1 ] - [r_{21}, L_2]
\label{rma}
\eeq

It must immediately be remarked that such an object is by no means unique. 
Moreover there is no one-to-one correspondance between a given dynamical
system and the Lie algebra in which its Lax representation is defined; a same
dynamical system may have several Lax representations and several $r$-matrix
structures.
\newline

{\bf 4- The classical Yang-Baxter equation}

The Poisson bracket structure (\ref{rma}) obeys a Jacobi identity which implies
an algebraic constraint for the $r$-matrix. Since $r$ depends a priori on
the dynamical variables this constraint takes a complicated form :

\beq
 [ L_1 ,[r_{12},r_{23}] + [r_{12}, r_{13}] + [r_{32},r_{13}]
+ \{ L_2, r_{13} \} - \{ L_3, r_{12} \} ] + \mbox{cycl. \ perm.} = 0
\label{yb2}
\eeq

Relevant particular cases of this very general identity are obtained when:

a) $r$ is independant of $x,q$. One is then lead to the general non-dynamical
Yang Baxter equation \cite{Skl1,JMM1,BV1,AT1}:

\beq
[r_{12},r_{23}] + [r_{12}, r_{13}] + [r_{32},r_{13}] = 0
\label{yb3}
\eeq

b) If furthermore $r$ is antisymmetric under permutation of the two copies
of the algebra ${\cal G}$ one obtains the better known and much studied
\cite{Skl1,BZ} form:

\beq
[r_{12},r_{23}] + [r_{12}, r_{13}] + [r_{13},r_{23}] = 0
\label{yb4}
\eeq

c) If on the contrary $r$ is dynamical, the supplementary 
terms in (\ref{yb2}) may take a completely algebraic form due to the specific
structure of the Lax operator. For instance the Calogero-Moser models 
lead to a self-contained algebraic equation for a classical $r$ matrix
depending only on position-type canonical variables $x_{\nu}$:

\beq
[r_{12},r_{13}]+[r_{12},r_{23}]+[r_{13},r_{23}]
-\sum_{\nu} h_{\nu}^{(1)} \frac{\partial}{\partial x_{\nu}} r_{23}
+\sum_{\nu} h_{\nu}^{(2)} \frac{\partial}{\partial x_{\nu}} r_{13}
-\sum_{\nu} h_{\nu}^{(3)} \frac{\partial}{\partial x_{\nu}} r_{12}=0
\label{yb5}
\eeq
where the set $  \{ h_{\nu}^{(i)} \}$ is a choice of basis for the Cartan
algebra of ${\cal G}$ acting on the representation space $i$. This equation
was first derived in \cite{FeWie} and a classification scheme of its solutions
was proposed in \cite{Schiff}, closely connected to the general algebraic 
scheme in \cite{BZ}.

Note by the way that there exists canonical examples of dynamical $r$ matrices,
obtained by using the well-known higher Poisson bracket construction for
any integrable dynamical system \cite{Skl2,Magri} starting from a constant
$r$-matrix. In particular the quadratic Sklyanin
bracket, where the Poisson structure of a Lax matrix becomes 
$\{ L_1 , L_2 \} = [r_{12}, L_1 \otimes L_2]$, may also be described as
a linear structure with a dynamical $r$ matrix $R_{12} \equiv r_{12} L_2$.
We shall see on the example of the Ruijsenaars-Schneider system that
dynamical ``linear'' $r$-matrix structures may also give
rise to dynamical ``quadratic'' structures. The {\it initial}
dynamical $r$-matrices themselves, however, are themselves {\it not} of the 
Sklyanin type.

\section{Hamiltonian reduction and $r$-matrices}

\subsection{General hamiltonian reduction}

We begin by recalling some well-known facts concerning the Hamiltonian 
reduction
of dynamical systems whose phase space is a cotangent bundle \cite{AM78}. Let 
$M$ be a
manifold and $N=T^* M$ its cotangent bundle. $N$ is equipped with the canonical
1-form $\alpha$ whose value at the point $p \in T^* M$ is $\pi^* p$ where 
$\pi$ is
the projection of $N$ on $M$.
If a Lie group $G$ acts on $M$, each
element $X \in \cal G$ (the Lie algebra of $G$) generates a vector field on $M$
that we shall denote $X.m$ at the point $m \in M$.
It lifts to a vector field on $N$ leaving $\alpha$
invariant. We shall also denote $X.p \in T_p(N)$ the value at $p \in N$ of this
vector field, so that the Lie derivative ${\cal L}_{X.p} \, \alpha$ of the 
canonical 1-form vanishes.  

$N$ is a symplectic manifold equipped with the
canonical 2-form $\omega = - {\rm d}\alpha $
To any function (or Hamiltonian) ${ \rm H}$
on $N$ we associate a vector field $X_{ \rm H}$ such that ${\rm d} 
{ \rm H} = i_{X_{ \rm H}}\, \omega $ 
and conversely since $\omega$ is non--degenerate.

The Hamiltonian associated to the above vector field $X.p,\; X \in \cal G$
reads:

\begin{equation} 
{ \rm H}_X(p)=i_{X.p}\, \alpha = \alpha\,(X.p) \label{moment}
\end{equation}

For any two functions $F$, $G$ on $N$  one defines the Poisson bracket $\{F,G\}$
as a function on $N$ by:
\begin{equation}
\{F,G\}=\omega\,(X_F,X_G)
\end{equation}
The Poisson bracket of the Hamiltonians associated to the group action has a
simple expression. 
In fact the group action is Poissonnian, i.e.:
\begin{equation}
\{{ \rm H}_X,{ \rm H}_Y\}={ \rm H}_{[X,Y]} \label{2.hom}
\end{equation}

Obviously the application $X\in {\cal G} \to { \rm H}_X(p),\; 
{\rm any}\; {\rm for any} p \in N$, is
a linear map from ${\cal G}$ to the scalars and so defines an element
${\cal P}(p)$ of ${\cal G}^*$ 
which is called the momentum at $p \in N$.

One then restricts oneself to the submanifold $N_\mu$ of $N$
with fixed momentum $\mu$ i.e. such that: $N_\mu ={\cal P}^{-1} (\mu) $

Due to  equation (\ref{moment})  and the invariance of $\alpha$
the action of the group $G$ on $N$ is 
transformed by $\cal P$ into the coadjoint action of $G$ on ${\cal G}^*$
\begin{equation}
{\cal P}(g.p)(X) = \alpha\,(g.g^{-1}Xg.p) =
{\rm Ad}^*_{g} \,{\cal P}(p) (X)
\end{equation}
where the coadjoint action on an element $\xi$ of ${\cal G}^*$ is
defined as:
$$ {\rm Ad}^*_g \, \xi \, (X) = \xi \, (g^{-1} X g) $$
The stabilizer $G_\mu$ of $\mu \in {\cal G}^*$ acts on $N_\mu$. The reduced 
phase space is precisely obtained by taking the quotient (assumed well-behaved):
\begin{equation}
{\cal F}_\mu = N_\mu / G_\mu
\end{equation}
It is known that this is a symplectic manifold. 

We then need to compute the Poisson bracket of functions 
on ${\cal F}_\mu$. These functions are conveniently described as $G_\mu$
invariant functions on $N_\mu$. To compute their Poisson
bracket we first extend them arbitrarily in the vicinity of $N_\mu$.
Two  extensions differ by a function vanishing on $N_\mu$. The difference
of the Hamiltonian vector fields of two such extensions
is controlled by the following:
\begin{lem}
Let $f$ be a function defined in a vicinity of $N_\mu$ and vanishing on $N_\mu$.
Then the Hamiltonian vector field $X_f$ associated to $f$ is tangent to the
orbit $G.p$ at any point $p \in N_\mu$. 
\end{lem}

As a consequence of this lemma we have a method to compute the reduced Poisson
bracket. We take two  functions defined on $N_\mu$ and invariant under $G_\mu$
and extend them arbitrarily. Then we compute their Hamiltonian vector fields on $N$
and project them on the tangent space to $N_\mu$ by adding a vector tangent
to the orbit $G.p$. These projections are independent of the extensions and the
reduced Poisson bracket is given by the value of the symplectic form on $N$ 
acting on them.
\begin{prop}
At each point $p \in N_\mu$ one can choose a vector $V_f.p \in {\cal G}.p$ such
that $X_f + V_f.p \in T_p(N_\mu)$ and $V_f.p$ is determined up to a vector
in ${\cal G}_\mu .p$.
\end{prop}

One finally gets the consistent general formula for reduced Poisson brackets:

\begin{prop}
The reduced Poisson bracket of two functions on ${\cal F}_\mu$ can be
computed using any extensions $f,\, g$ in the vicinity of $N_\mu$
according to:
\begin{equation}
\{f,g\}_{\rm reduced}=\{f,g\}+{1\over 2}\left(
(V_g.p).f -(V_f.p).g
\right) \label{bienreduit}
\end{equation}
This is equivalent to  the Dirac bracket. 
\end{prop}

\subsection{The case $N=T^*G$}

Let now $M=G$ be a Lie group; one uses the left translations to
identify $N=T^*G$ with $G\times {\cal G}^*$.
\begin{equation}
\omega \in T^*_g(G) \longrightarrow (g,\xi)\mbox{~~~where~~}
\omega=L^*_{g^{-1}}\xi
\end{equation}

The Poisson structure on $N=T^*G$ is easily seen to be:

\beq 
\{ \xi (X), \xi (Y) \} = - \xi ([X,Y]) \,\, ; \,\, 
\{ \xi (X) ,g \} = - g\;X \,\, ; \,\,
\{ g , g \} =0 \label{poi3}
\eeq

Geodesics on the group $G$ correspond to  left translations of 
1-parameter groups (the
tangent vector is transported parallel to itself), therefore
${d\over dt}(g^{-1}\dot{g}) =0$.
This is a Hamiltonian system whose Hamiltonian is:
${ \rm H} = {1 \over 2}\; (\xi,\xi)$.
where we have identified ${\cal G}^*$ and $\cal G$ through the invariant
Killing metric.

Here ${\rm H}$ is bi--invariant, so one can 
reduce this dynamical system using  Lie subgroups $H_L$ 
and $H_R$ of $G$ of Lie algebras ${\cal H}_L$
and ${\cal H}_R$, acting respectively on the left and on the right on $T^*G$
in order to obtain a non--trivial result.

Using the coordinates $(g,\xi)$ on $T^*G$  this action reads:
$$(\;(h_L,h_R),(g,\xi)\;)\to (h_L g h_R^{-1},{\rm Ad}^*_{h_R} \xi)$$
We have written this action as a left action on $T^*G$, in order to
apply the formalism developed in Section {\bf 3.1}.

The moments are:
\begin{equation}
{\cal P}^L(g,\xi ) = P_{{\cal H}^*}\, {\rm Ad}^*_{g}\, \xi \, ; \quad 
{\cal P}^R(g,\xi ) = - P_{{\cal H}^*}\, \xi \label{1.mom} \, ; \quad
{\cal P}=({\cal P}^L,{\cal P}^R)
\end{equation}
where we have introduced the projector on ${\cal H}^*$ of forms in
${\cal G}^*$ induced by the restriction of these forms to $\cal H$.

\subsection{The Calogero-Moser models}

We follow here the derivation
of \cite{OlPer}.
Let us consider an involutive automorphism $\sigma$ of a simple Lie group $G$ 
and
the subgroup ${ H}$ of its fixed points. Then $H$ acts on the right on $G$ defining a
principal fiber bundle of total space $G$ and base $G/H$, which is a global
symmetric space. Moreover $G$ acts on the left on $G/H$ and in particular so 
does
$H$ itself. We shall consider the situation described in Section {\bf 3.2} when
$H_L=H_R=H$. The Hamiltonian of the geodesic flow on $G/H$ is
invariant under the $H$ action allowing to construct the Hamiltonian reduction
which under suitable choices of the momentum leads to the Calogero--Moser 
models. 
As a matter of fact since the phase space of the Calogero model is non-compact
one has to start from a non-compact Lie group $G$ and quotient it by a
maximal compact subgroup $H$ so that the symmetric space $G/H$ is of the
non-compact type. 
The derivative of $\sigma$ at the unit element of $G$ is an involutive automorphism
of $\cal G$ also denoted $\sigma$. Let us consider its eigenspaces $\cal H$ and
$\cal K$  associated with the eigenvalues $+1$ and $-1$ respectively. Thus we
have a decomposition:
\begin{equation}
{\cal G}={\cal H} \oplus {\cal K} \label{ghk}
\end{equation}
in which $\cal H$ is the Lie algebra of $H$ which acts by inner automorphisms on
the vector space $\cal K$ ($h{\cal K}h^{-1}=\cal K$).

 Let $\cal A$ be a
maximal commuting set of elements of $\cal K$. It is called a Cartan
algebra of the symmetric space $G/H$. It is known that every element in
$\cal K$ is conjugated to an element in $\cal A$ by an element of $H$.
Moreover $\cal A$ can be extended to a maximal commutative subalgebra of
$\cal G$ by adding to it a suitably chosen 
abelian subalgebra $\cal B$ of $\cal H$. We shall
use the radicial decomposition of $\cal G$ under the abelian 
algebra  
${\cal A}$:
\begin{equation}
{\cal G}={\cal A}\bigoplus{\cal B}\bigoplus_{e_\alpha,\; 
\alpha \in \Phi}\,\Bbb{R} e_\alpha
\label{radix}
\end{equation}

These decompositions of $\cal G$ exponentiate to similar decompositions of $G$.
First $G=KH$ where $K=\exp ({\cal K})$. 
Then $A=\exp ({\cal A})$ is a maximal
totally geodesic flat submanifold of $G/H$ and any element of $K$ can be
written as $k=h Q h^{-1}$ with $Q \in A$ and $h \in H$. It  follows that any
element of $G$ can be written as $g=h_1 Q h_2$ with $h_1,h_2 \in H$. 

Of course this decomposition is  non unique. This non--uniqueness
is described  in the following:
\begin{prop} \label{unicite}
If $g=h_1Qh_2=h'_1Q'h'_2$ we have: $h'_1=h_1d^{-1}h_0^{-1}$,
$h'_2=h_0dh_2$ and $Q'=h_0Qh_0^{-1}$ where $d\in \exp({\cal B})=B$
and $h_0\in H$ is a representative of an element of the Weyl group
of the symmetric space. So if we fix $Q=\exp(q)$ such that $q$ be in a 
fundamental Weyl chamber, the only ambiguity resides in the element $d \in B$.
\end{prop}

The reduction to Calogero-Moser models is then obtained by an adequate 
choice of the momentum
$\mu=(\mu^L,\mu^R)$ such that ${\cal P}=\mu$. We take $\mu^R=0$
so that the isotropy group of the right component is $H_R$ itself.

The choice of the moment $\mu^L$ is of course of crucial
importance. It must be fixed so that:

\begin{itemize}
\item  its isotropy group $H_\mu$ is a maximal proper Lie subgroup of $H$,
so that the phase space of the reduced system be of minimal dimension but 
non trivial.
\item In order to ensure the unicity of the decomposition introduced
in the Proposition~(\ref{unicite})
on $N_\mu$ we need:
\begin{equation}
{\cal H}_\mu \cap {\cal B}=\{ 0\}  \label{2.cap}
\end{equation}

We choose a complementary maximal isotropic subspace
$\cal C$ so that 
\begin{equation}
{\cal H}={\cal H}_\mu \oplus {\cal B} \oplus \cal C \label{hbc}
\end{equation}
and $\chi$ is a non--degenerate skew--symmetric bilinear form on
${\cal B} \oplus \cal C$, hence ${\rm dim }\,{\cal B}={\rm dim }\,{\cal C}$.
Notice that $\cal C$ is defined up to a symplectic transformation
preserving $\cal B$.
\item The reduced phase space ${\cal F}_\mu$ has dimension $2\,{\rm dim}\,
{\cal A}$
\end{itemize}

We now construct a section
$\cal S$ of the bundle $N_\mu$ over ${\cal F}_\mu$ so that one can write:
\begin{equation}
N_\mu=H_\mu {\cal S}H \label{ndec}
\end{equation}
To construct this section we  take a point $Q$ in $A$ and
an $L \in {\cal G}^*$ such that the point $(Q,L)$ is in $N_\mu$. 
In this subsection we shall for convenience identify $\cal G$ and
${\cal G}^*$ under the Killing form assuming that $G$ is semi--simple.
Moreover since the automorphism $\sigma$ preserves the Killing form,
$\cal H$ and $\cal K$ are orthogonal, and $P_{{\cal H}^*}$ reduces to
the orthogonal projection on $\cal H$.
Since $\mu^R=0$ we have $L\in {\cal K}$ and one can write:

\begin{equation}
L=p+\sum_{e_\alpha,\;\alpha \in \Phi'} l_\alpha\, (e_\alpha - \sigma(e_\alpha)) \label{2.L}
\end{equation}
where $p \in \cal A$.
From equation (\ref{1.mom}) one gets:
$$\mu^L= P_{{\cal H}} \left(
p + \sum_\alpha l_\alpha\, (Q e_\alpha Q^{-1} - Q \sigma(e_\alpha) Q^{-1})
\right)$$
Since $Q = \exp(q),\; q\in \cal A$ we have  $Qe_\alpha Q^{-1}=\exp(\alpha(q))
 e_\alpha$ and similarly $Q\sigma(e_\alpha) Q^{-1}=\exp(-\alpha(q)) 
\sigma(e_\alpha)$ . 

Then the above equation becomes:
\begin{equation}
\mu^L=\sum_{\alpha} l_\alpha \sinh \alpha(q)\, (e_\alpha +\sigma(e_\alpha))
\label{2.mu}
\end{equation}
One can choose the momentum of the form:
$\mu^L = \sum_\alpha g_\alpha (e_\alpha + \sigma(e_\alpha))$ namely $\mu^L$
has no component in $\cal B$, where the $g_\alpha$ are such
that $H_\mu$ is of maximal dimension (we shall see that it essentially
fixes them, and obviously if $g_\alpha \neq 0$ for any $\alpha$
equation~(\ref{2.cap}) is automatically satisfied) and we have shown the:
\begin{prop} \label{laxpair}
The couples $(Q,L)$ with $Q =\exp(q)$ and
$$L=p+\sum_{\alpha }\,{g_\alpha \over \sinh \alpha(q)}\, 
(e_\alpha-\sigma(e_\alpha))$$
with $p,q \in \cal A$ form a submanifold in $N_\mu$ of dimension $2\,{\rm dim}\,
{\cal A}$.
\end{prop}

Notice that $L$ is just the Lax operator of the Calogero model and that
the section $\cal S$ depends of $2\,{\rm dim}\,{\cal A}$ parameters in an
immersive way. 
Hence one can identify $N_\mu$ with the set of orbits of $\cal S$
under $H_\mu \times H$ i.e. the set of points
$(g=h_1Qh_2,\xi=h_2^{-1}Lh_2)$ with $h_1\in H_\mu$ and $h_2\in H$
{\em uniquely defined} due to condition~(\ref{2.cap}). The variables
$p$ and $q$ appearing in $Q$ and $L$ are the dynamical variables of
the Calogero model and form a pair of canonically conjugate variables.

We then  compute the Poisson bracket of the functions on ${\cal F}_\mu$
whose expressions on the section $\cal S$ are $L(X)$ and $L(Y)$ for
$X,Y \in \cal K$. These functions have uniquely defined $H_\mu \times H$
invariant extensions to $N_\mu$ given respectively by:
$$F_X(g,\xi)=<\xi,h_2^{-1}Xh_2>,\;F_Y(g,\xi)=<\xi,h_2^{-1}Yh_2>
{\rm ~where~} g=h_1Qh_2$$
Notice that $h_2$ is a well--defined function of $g$ in $N_\mu$ due to
condition~(\ref{2.cap}). According to the prescription given in the
section {\bf 3.1} we choose extensions of these functions in the vicinity
of $N_\mu$. We {\em define} these extensions at the point $p=(g,\xi) \in T^*G$
by the {\em same} formulae in which $h_2$ is chosen to be a function
depending {\em only} on $g$ and reducing to the above--defined $h_2$
when $p \in N_\mu$.
Because of the non--uniqueness of the decomposition $g=h_1Qh_2$
outside  $N_\mu$ one cannot assert that the functions $F_X,\,F_Y$ are
invariant under the action of $H\times H$ and we must appeal to the
general procedure to compute the reduced Poisson brackets.

The complete derivation with all its technical subtleties can be found in
\cite{AvBabTa}. The final result gives the general $r$-matrix for trigonometric
Calogero-Moser models in so-called dual form \cite{Skl1}:

\begin{th}\label{Rgene}
There exists a linear mapping $R:\;{\cal K}\to {\cal H}$ such that:
\begin{equation}
\{ \, L(X)\, ,\, L(Y)\, \}_{\rm reduced}=L\, (\,[X,RY]+[RX,Y]\,) \label{matricer}
\end{equation}
and $R$ is given by:
\begin{equation}
R\,(X)= \nabla_g h_2\,(X) + {1 \over 2} D_Q\,( V_X)
\label{rformel}
\end{equation}
\end{th}
where:
\begin{prop}\label{nablah2}
On the section $\cal S$ with $Q=\exp(q) \in A$ we have:\hfil\break
For $X \in \cal K$ i.e. $X=X_0+\sum X_\alpha (e_\alpha-\sigma e_\alpha),
\; X_0 \in \cal A$
\begin{equation}
\nabla_g h_2  ( X)=-h_0(X)+\sum_\alpha X_\alpha \coth (\alpha(q))
(e_\alpha+\sigma e_\alpha) \label{dh2}
\end{equation}
Here $h_0(X)$ is a linear function from $\cal G$ to $\cal B$ which 
is fixed by the condition:
\begin{equation}
X_L \equiv h_0(X)-\sum_{\alpha }\,{X_\alpha \over \sinh \alpha(q)}
\,(e_\alpha+\sigma e_\alpha) \in {\cal H}_\mu \oplus \cal C
\label{mmain}
\end{equation}
\end{prop}

\subsection{Two examples}

To illustrate the power of this method we now give two examples of $r$-matrices.
The $A_n$ case had already been treated in \cite{AvTaSkBraSu} and serves as a
check on the validity of the derivation. The case of $SU(n,n)$ Calogero-Moser
model proved to be too intricate for a direct computation; however this method
immediately gives its $r$-matrix.
\newline

{\bf The standard Calogero-Moser model} ($SL(n)$)
is obtained by starting from
the non compact group $G=SL(n,\Bbb{C})$ and its maximal compact
subgroup $H=SU(n)$ as first shown by~\cite{OlPer}. We
choose the momentum $\mu_L$ as described in Section~(3.1) so that
the isotropy group $H_\mu$ be a maximal proper Lie subgroup of $H$.
Obviously one can take $\mu_L$ of the form:
\begin{equation}
\mu_L= i\,(vv^+ -1) \label{musln}
\end{equation}
where $v$ is a vector in $\Bbb{C}^n$ such that $v^+ v=1$, hence $\mu_L$
is a traceless antihermitian matrix. Then $g \mu_L g^{-1} =\mu_L$ if
and only if $g v= c v$ where $c$ is a complex number of modulus 1. Hence
$H_\mu = S(U(n-1)\times U(1))$ which has the above--stated property.

In this case the automorphism $\sigma$ is given by $\sigma\,(g)=
(g^+)^{-1}$ (notice that we consider only the real Lie group structure),
$B$ is the group of  diagonal matrices of determinant 1 with 
pure phases on the diagonal and $A$ is the group of real diagonal
matrices with determinant 1. The property~(\ref{2.cap})
is then satisfied as soon as the vector $v$ has no zero component.
As a matter of fact, $v$ is further constrained by $\mu_L$
being a value of the
moment map. Considering equation~(\ref{2.mu}) we see that $\mu_L$ has
no diagonal element, which implies that all the components of $v$
are pure phases $v_j=\exp (i\theta_j)$. 
These extra phases which will appear in the Lax matrix
can however be conjugated out by the adjoint action of a constant matrix
${\rm diag}\,(\exp (i\theta_j))$ hence we shall from now on set $v_j=1$
for all $j$. This is the solution first considered by Olshanetskii
and Perelomov.

The Lax matrix $L$ is then given by Proposition~(\ref{laxpair})
and therefore
$$L=p+\sum_{k<l}\, {1\over \sinh(q_k-q_l)} (i E_{kl} - i E_{lk})$$

The $r$--matrix can now be deduced straightforwardly from 
Proposition~(\ref{nablah2}), after reconverting the dual form where $R$ is an
endomorphism of the Lie algebra into
the more usual direct form where $R$ lives in the tensor product of
the Lie algebra by itself. One ends up with:

$$R_{12}=\sum_{k\neq l}\, \coth (q_k-q_l) E_{kl}\otimes E_{lk}
+ {1 \over 2} \sum_{k \neq l}\, {1 \over \sinh (q_k-q_l)}
(E_{kk}-{1 \over n}\,{\bf 1})\otimes (E_{kl}-E_{lk})$$
This gives back the already known $r$--matrix of the Calogero model for the
potential $1/\sinh(x)$, and the other potentials $1/\sin(x)$ and
$1/x$ have similar $r$--matrices obtained by analytic continuation.
\newline

{\bf The $SU(n,n)$ Calogero model} is obtained by starting from the 
non-compact group $G=SU(n,n)$. This is the subgroup of $SL(2n,{\Bbb C})$ 
which leaves invariant the sesquilinear quadratic form defined by
\begin{eqnarray}
Q((u_1,v_1),(u_2,v_2)) = \pmatrix{u_1^+ & v_1^+} J
\pmatrix{u_2 \cr v_2}= u_1^+ v_2 + v_1^+ u_2
\label{quadra}
\end{eqnarray}
where $u_i,v_i$ are vectors in ${\Bbb C}^n$ and $J$ is the matrix
\begin{eqnarray} 
J= \pmatrix{0 & {\bf 1} \cr {\bf 1} & 0}.
\nonumber
\end{eqnarray}
The Lie algebra of $SU(n,n)$ therefore
consists of block matrices
\begin{eqnarray}
{\cal G}= \{ \pmatrix{a & b \cr c & d } \vert a=-d^+ ,\; {\rm Tr}\; (a+d)=0,
\; b^+ = -b,\; c^+ = -c \}
\label{sunn}
\end{eqnarray}
where $a,b,c,d$ are $n \times n$ complex matrices.

We consider again the automorphism $\sigma: \sigma(g)=(g^+)^{-1}$, which can be
consistently restricted to $SU(n,n)$. Its fixed points at the Lie algebra level
consist of block matrices
\begin{eqnarray} 
{\cal H}=\{ \pmatrix{a & c \cr c & a} \vert a^+ =-a, \; { \rm Tr } \,  (a)=0,\; c^+ =-c \}
\label{sunun}
\end{eqnarray}
This Lie algebra is isomorphic to the Lie algebra of $S(U(n) \times U(n))$, the two
$u(n)$'s being realized respectively by $a+c$ and $a-c$.

The subalgebra ${\cal B}$ consists of matrices of the form (\ref{sunun}) with
$c=0$, and $a$ is a diagonal matrix of zero trace and purely imaginary
coefficients. The Abelian subalgebra ${\cal A}$ consists of matrices of the form
(\ref{sunn}) with $b=c=0$ and $a=-d$ being a real diagonal matrix.

To perform the reduction, we choose as above $\mu^R=0$ and 
\begin{eqnarray} 
\mu^L = i(v v^+ - {\bf 1})+ i\gamma J
\label{momentsu}
\end{eqnarray}
The vector $v$ has again $2n$ components all equal to 1.

Notice that in equation (\ref{momentsu})
the parameter $\gamma$ is an arbitrary real number. This will
lead to existence of a second coupling constant in the
corresponding Calogero model.

Then, from proposition (\ref{laxpair}), the Lax matrix is found to be:

\begin{eqnarray} 
L&=& p +\sum_{i<j}{1\over \sinh (q_i -q_j)}(1-\sigma)
(iE_{ij}+iE_{j+n,i+n}) \nonumber \\
&&+  \sum_{i<j}{1\over \sinh (q_i +q_j)}(1-\sigma)
(iE_{i,j+n}+iE_{j,i+n}) \nonumber \\
&&+ (\gamma +1)\sum_i{1\over \sinh (2 q_i )} (1-\sigma)(iE_{i,i+n})
\nonumber
\end{eqnarray}
where $p$ is a generic element of ${\cal A}$ of the form ${\rm
diag}\, p_i,-{\rm diag}\, p_i$.

The $r$-matrix is then computed straightforwardly:
\begin{eqnarray} 
R_{12}&=&{1\over 2} \sum_{k\neq l} \coth (q_k-q_l) (E_{kl}+E_{k+n,l+n})
\otimes (E_{lk}-E_{l+n,k+n}) \nonumber \\
&+&{1\over 2} \sum_{k,l} \coth (q_k +q_l) (E_{k,l+n}+E_{k+n,l})
\otimes (E_{l+n,k}-E_{l,k+n}) \nonumber \\
&+& {1\over 2} \sum_{k\neq l}{1\over \sinh (q_k -q_l)}
(E_{kk}+E_{k+n,k+n}-{1\over n} {\bf 1})\otimes (E_{kl}-E_{k+n,l+n})
\nonumber \\
&+&{1\over 2} \sum_{k,l}{1\over \sinh (q_k+q_l)}
(E_{kk}+E_{k+n,k+n}-{1\over n} {\bf 1})\otimes (E_{k,l+n}-E_{k+n,l})
\nonumber 
\end{eqnarray}
 
These dynamical $r$-matrices depend only on the dynamical variable $q$. The 
different approach advocated in \cite{Bra2} leads to $r$-matrices depending
on both $p$ and $q$ variables, but on a smaller set of algebra generators.

\section{The dynamical $r$ matrices of Calogero and Ruijsenaars models}

Dynamical $r$-matrices have been derived for the Calogero-Moser and 
(relativistic) Ruijsenaars Schneider models using either the technique described
here or a direct method starting from an ansatz of the
same form. We will now describe the results achieved in this way for $A_n$
models, and indicate interesting and sometimes deep connections between these 
various $r$-matrices.

Let us start with {\bf Calogero-Moser models}. The rational and 
trigonometric matrices
were described in the previous section.. The elliptic case was solved by
Sklyanin and by Braden et al. \cite {AvTaSkBraSu}. The Lax matrix reads:
\beq
  L(\la)  =  \sumi p_i \ e_{ii} + \sumij l(q_{ij},\la)  \ e_{ij}
\label{eqq}
\eeq
Here one has set:
\beq
  l(x,\la) = - \frac{\si(x+\la)}{\si(x)\ \si(\la)}, \ \ \ \ \ 
  V(x) = \wp(x)
\eeq
where $\si$ and $\wp$ are Weierstrass elliptic functions. The classical 
$r$-matrix reads:

\bea
   r_{12}(\la,\mu) & = &
     \sumij l(q_{ij},\la-\mu) \  e_{ij} \otimes e_{ji}                    
   + \demi \ \sumij l(q_{ij},\mu) \ (e_{ii}+e_{jj}) \otimes e_{ij}
         \nonumber \\
  & &  - [\zeta(\la-\mu) + \zeta(\mu)] \sumi e_{ii} \otimes e_{ii}.
\label{ells} 
\eea

Note that a spectral parameter is now present in $L$ and $r$. This particular
formulation of the spinless elliptic case is due to Krichever \cite{Kri}. The
other known Lax formulation due to Olshanetskii and Perelomov \cite{OlPer} has
no spectral parameter but requires a $p$ and $q$ dependance in the $r$-matrix,
which was only recently given in \cite{Bra2} and has a totally different
algebraic form.
\newline

{\bf The spin Calogero-Moser models} were introduced in \cite {GibHe}. 
The Lax operator for the elliptic case reads:

\beq
  L(\la)  =  \sumi p_i \ e_{ii} + \sumij l(q_{ij},\la) \ f_{ij} \ e_{ij}
\label{eq2}
\eeq

where $f_{ij}$ are spinlike variables with the Kirillov-Poisson bracket 
structure:
\beq
\{ f_{ij}, f_{kl} \}  =  \demi \ ( \delta_{il} \ f_{jk} +
                                 \delta_{ki} \ f_{lj} +
                                 \delta_{jk} \ f_{il} +
                                 \delta_{lj} \ f_{ki}  ) .
\eeq

One then needs to introduce a parametrization of $f_{ij}$ so as to be on a coadjoint
orbit of $SU(N)$:
Introducing vectors
\bean  
  (\xi_i)_{i=1 \cdots N} & \ \mbox{with} & \xi_i=(\xi_i^a)_{a=1 \cdots r } \\
  (\eta_i)_{i=1 \cdots N} & \ \mbox{with} & \eta_i=(\eta_i^a)_{a=1 \cdots r }
\eean 
with the Poisson brackets
\beq
  \{ \xi_i^a ,\xi_j ^b \} = 0 , \ \ \ \ \ \{ \eta_i^a , \eta_j^b \} = 0 , 
  \ \ \ \ \ \{ \xi_i^a ,\eta_j^b \} = - \delta_{ij} \ \delta_{ab} ,
\eeq
we parametrize $f_{ij}$ as follows:
\beq
  \label{sm}
  f_{ij}= \langle  \xi_i | \eta_j \rangle = \sum_{a=1}^r \xi_i^a \eta_j^a.
\eeq
The phase space now becomes a true symplectic manifold.

The Hamiltonian takes the form
\beq
  \label{Ham}
  H = \demi \sumi p_i^2 - \demi \sumij f_{ij} \ f_{ji} \ V(q_{ij}) ,
  \ \ \ \ \  q_{ij}=q_i-q_j  
\eeq
with the Weierstrass function as elliptic potential, as in the spinless case.

The classical $r$-matrix then reads:

\bea
   r_{12}(\la,\mu) & = & \demi \ \sumij l(q_{ij},\la-\mu)e_{ij} 
                         \otimes e_{ji} 
                    + \demi \ \sumij l(q_{ij},\la+\mu)e_{ij} 
                         \otimes e_{ij} \nonumber \\
                    & - & \!  \demi \ [ \zeta(\la+\mu)+\zeta(\la-\mu) ] 
                         \sumi  e_{ii} \otimes e_{ii}. 
\eea

Trigonometric and rational cases can be derived from the elliptic case by 
taking suitable limits \cite{BiAvBa}.
The spinless case can also be derived from the spin case by taking $r=1$ and
introducing a further Hamiltonian reduction by the action of $U(1)$ as
a phase on the vectors $\xi_i, \eta_i$. The supplementary terms in the spinless
$r$ matrix arise from the conjugation of the Lax matrix required to bring it
in canonical shape (\ref{eqq})  after elimination of the vector-like degrees 
of freedom. Let
us finally remark that this $r$-matrix structure has yielded a number of
important developments: exact classical Yangian symmetry \cite{BiAvBa} 
(a quantum version of it had been found beforehand, using heavy direct 
algebraic computations \cite{BeGauHaPa}); quantum
version of the dynamical $r$-matrix using the shifted version of the
quantum Yang-Baxter equation described in \cite{MaNi,GeNe,Fel1}.
\newline

{\bf The spinless relativistic RS models} are described by the hamiltonian:

\beq
H = mc^2 \sum_{j=1}^N \left( {\rm cosh} \, \theta_j \right) \prod_{k\not= j} \,
f (q_k - q_j)
\eeq
where
\bea
f(q) & = & \left(1 + {g^2\over q^2} \right)^{1/2}\,\, ({\rm rational})\nonumber
\\
f(q) & = & \left( 1 + {\alpha^2\over {\rm sinh}^2 {\nu q\over
2}}\right)^{1/2}\,\,{\rm (hyperbolic)}\nonumber\\
&&\nonumber\\
f (q) & = & \left( \lambda + \nu {\cal P} (q) \right) \,\, {\rm ( elliptic)},
\,\,{\cal P} = {\rm Weierstrass\,\, function} \label{poten}
\eea

Here the canonical variables are a set of rapidities 
$\{ \theta_i, i=1\cdots N\}$
and conjugate positions $q_i$ such that $\{ \theta_i, q_j\} = \delta_{ij}$.

The dynamical system admits a Lax representation with the Lax 
operator:

\bea
&L =\sum_{j,k=1}^N \, L_{jk} \, e_{jk} \nonumber\\
&L_{jk} = \exp {\beta}  \theta_j  \cdot C_{jk}
\left( q_j - q_k \right) \cdot
\left( \displaystyle\prod_{m\not= j} \, f \left( q_j - q_m \right)
\prod_{l\not= k} \, f \left( q_l - q_k \right)\right)^{1/2}
\eea

where $\{ e_{jk} \}$ is the usual basis for $N \times N$ matrices; $f$ was
given in (\ref{poten}) and
\bea
C_{jk} (q) & = & {\gamma\over\gamma + iq} \,\,\,\,{\rm (rational)}\nonumber\\
C_{jk} (q) & = & \left( {\rm cosh} {\nu\over 2} q + ia \, {\rm sinh} {\nu\over
2} q \right)^{-1} \,\,\,\,{\rm (trigonometric)}\\
C_{jk} (q) & = & \frac{\Phi(q + \gamma, \lambda)}{\Phi(\gamma, \lambda)}
\,\, \,\,{\rm (elliptic)} \label{factor}
\eea

Again the function $\Phi$ is defined as:

\beq
\Phi (x, \lambda) \equiv \frac{\sigma(x+\lambda)}{\sigma(x)\sigma(\lambda)}
\eeq
where $\sigma$ is the Weierstrass function.

The elliptic $r$-matrix structure is better written as a quadratic expression
in terms of the Lax operator \cite{Su}:

\bea
\{ L_1(\lambda) , L_2 (\mu)\} = && ( L_1(\lambda) \otimes L_2 (\mu))
a_1(\lambda, \mu) - a_2(\lambda, \mu) ( L_1(\lambda) \otimes L_2 (\mu))
\nonumber\\
&&+ ({\bf 1} \otimes L_2(\mu)) s_1(\lambda, \mu) (L_1(\lambda) \otimes {\bf 1})
\nonumber\\
&&- (L_1(\lambda) \otimes {\bf 1})  s_2(\lambda, \mu) 
({\bf 1} \otimes L_2(\mu)) \label{quad}
\eea

Here one defines:

\bea
a_1(\lambda, \mu) = a(\lambda, \mu) +w \,\, &,& \,\, 
s_1(\lambda, \mu) = s(\lambda, \mu) -w \nonumber\\
a_2(\lambda, \mu) = a(\lambda, \mu) +s(\lambda ) -s^{\ast} (\mu) &-& w \,\, , \,\,
s_2(\lambda, \mu) = s^{\ast} (\mu) +w \label{quad2}
\eea

The matrices $a$ and $s$ are obtained from the $r$-matrix of the elliptic
Calogero-Moser model given in (\ref{ells}) as $r(\lambda, \mu) \equiv 
a(\lambda, \mu) + s(\lambda)$ where $a$ is the skew-symmetric matrix:
\beq
a(\lambda, \mu) = -\zeta(\lambda - \mu)\Sum_{k=1}^{N} E_{kk} \otimes E_{kk}
-\sum_{k \neq j} \Phi(q_j -q_k, \lambda - \mu ) E_{jk} \otimes E_{kj}
\label{decomp1}
\eeq
and $s, s^{\ast}$ are  non-skew-symmetric matrices independent of the second
spectral parameter:

\bea
s(\lambda, \mu) &=& \zeta(\lambda)\sum_{k=1}^{N} E_{kk} \otimes E_{kk}
+ \Sum_{k \neq j} \Phi(q_j -q_k, \lambda ) E_{jk} \otimes E_{kk} \nonumber\\
s^{\ast}(\lambda, \mu) &=& \zeta(\lambda)\sum_{k=1}^{N} E_{kk} \otimes E_{kk}
+ \Sum_{k \neq j} \Phi(q_j -q_k, \lambda ) E_{kk} \otimes E_{jk} 
\label{decomp2}
\eea

and finally $w$ is a supplementary matrix, independent of the spectral 
parameters:

\beq
w = \Sum_{k \neq j} \zeta (\q_k -q_j) E_{kk} \otimes E_{jj}
\label{decomp3}
\eeq

This $r$-matrix structure is a Sklyanin-type bracket (although 
realized in the more generic case of an initially dynamical $r$-matrix) obtained
from the Calogero-Moser $r$-matrix structure viewed as a linear bracket. This
can be interpreted from the fact that RS models are obtained not only as
hamiltonian reductions from current algebras on elliptic curves but also
alternatively \cite{ACF1} as hamiltonian 
reductions from Heisenberg double of Lie groups \cite{Drin}. In this case the 
initial Poisson structure on the large phase space is itself a quadratic bracket
instead of the canonical initial linear (Kirillov) bracket which is
the natural structure on the cotangent bundle of a Lie group. This relation
is maintained throughout the hamiltonian reduction procedure and the final
$r$-matrix structures are essentially connected in the same way.

The previously obtained $r$-matrices \cite{AvRol} can be obtained from this
one by sending one period of the elliptic functions to infinity and suitably
conjugating the Lax pair in such a way as to get a completely symmetric
expression in terms of the momenta $\theta_i$. On the other hand the
$r$-matrix found in \cite{BaBe} cannot be easily inserted in this scheme. In 
fact it correspond to a very specific value of the parameters where the
Lax matrix becomes completely symmetric and tne $r$-matrix may only then
take this very special form.

The classical $r$-matrices admit a quantization scheme on the same lines as
the Calogero-Moser case \cite{ABB2}. It was developed in \cite{ACF2}.

Finally a word about {\bf the spin RS dynamical system}. They were introduced in
\cite{KriZa}. It is not clear at this time how to define a consistent
hamiltonian structure in the most general case although the rational case
was solved recently \cite{Russes}. The task is indeed easier here since there
exists a duality symmetry \cite{Rui2} connecting the rational RS model
to the trigonometric CM model for which the spin model is well known. Let
us finally mention that there exists a general scheme to obtain Hamiltonian 
structures from Lax representations using the tools of algebraic geometry
\cite{KriPho} (see also Pr. Krichever's contribution to this colloquium) 
and this scheme now appears to be the most promising way to get these
elusive Hamiltonian structures.

\vspace{1cm}

{\bf Acknowledgements}

The works presented here were done in collaboration with O. Babelon, E. Billey,
G. Rollet and M. Talon at the LPTHE Paris VI (URA CNRS 280). I have also drawn
upon works by H.W. Braden, Yu. B. Suris and E.K. Sklyanin. I wish to thank 
A. Antonov, H.W. Braden and J.M. Maillet for clarifications of particular 
points in this presentation; S.N.M. Ruijsenaars for suggesting to write a 
summary of the current situation as Section 4; and J.F. Van Diejen and 
L. Vinet for their kind invitation to attend this workshop.

\end{document}